\documentclass[12pt]{article}
\usepackage{mathrsfs}
\usepackage{amscd}
\usepackage{amsmath,amsfonts,amssymb,amscd}
\usepackage{indentfirst,graphics,epsfig,psfrag}
\input{epsf}
\usepackage{ifpdf}
\usepackage{enumerate}
\usepackage{appendix}
\usepackage{enumerate}

\usepackage{lineno}

\makeatletter \@addtoreset{figure}{section} \makeatother
\makeatletter
\long\def\@makecaption#1#2{%
   \vskip 10\p@
   \setbox\@tempboxa\hbox{{#1}\ \ #2}%
   \ifdim \wd\@tempboxa >\hsize
       {#1}\ \ #2\par
   \else
       \hbox to\hsize{\hfil\box\@tempboxa\hfil}%
   \fi}
\makeatother

\newtheorem{thm}{Theorem}[section]

\newtheorem{pro}[thm]{Proposition}

\newcommand{\qed}{{\hfill\rule{3pt}{7pt}}}

\def\qed{\hfill \rule{4pt}{7pt}}

\begin{document}
\title{\bf Note on the Rainbow $k$-Connectivity
of\\ Regular Complete Bipartite Graphs\footnote{Supported by NSFC,
PCSIRT and the ``973" program. }}
\author{
\small  Xueliang Li,  Yuefang Sun\\
\small Center for Combinatorics and LPMC-TJKLC\\
\small Nankai University, Tianjin 300071, P.R. China\\
\small E-mail: lxl@nankai.edu.cn; syf@cfc.nankai.edu.cn }

\date{}
\maketitle
\begin{abstract}
A path in an edge-colored graph $G$, where adjacent edges may be
colored the same, is called a rainbow path if no two edges of the
path are colored the same. For a $\kappa$-connected graph $G$ and an
integer $k$ with $1\leq k\leq \kappa$, the rainbow $k$-connectivity
$rc_k(G)$ of $G$ is defined as the minimum integer $j$ for which
there exists a $j$-edge-coloring of $G$ such that any two distinct
vertices of $G$ are connected by $k$ internally disjoint rainbow
paths. Denote by $K_{r,r}$ an $r$-regular complete bipartite graph.
Chartrand et al. in ``G. Chartrand, G.L. Johns, K.A. McKeon, P.
Zhang, The rainbow connectivity of a graph, Networks 54(2009),
75-81" left an open question of determining an integer $g(k)$ for
which the rainbow $k$-connectivity of $K_{r,r}$ is $3$ for every
integer $r\geq g(k)$. This short note is to solve this question by
showing that $rc_k(K_{r,r})=3$ for every integer $r\geq
2k\lceil\frac{k}{2}\rceil$, where $k\geq 2$ is a positive integer.\\[2mm]
\end{abstract}

\section{Introduction}

All graphs considered in this paper are simple, finite and
undirected. Let $G$ be a nontrivial connected graph with an edge
coloring $c: E(G)\rightarrow \{1,2,\cdots,k\}$, $k\in \mathbb{N}$,
where adjacent edges may be colored the same. A path of $G$ is
called $rainbow$ if no two edges of it are colored the same. A
well-known result shows that in every $\kappa$-connected graph $G$
with $\kappa \geq 1$, there are $k$ internally disjoint $u-v$ paths
connecting any two distinct vertices $u$ and $v$ for every integer
$k$ with $1\leq k\leq \kappa$. Chartrand et al. \cite{Chartrand 2}
defined the $rainbow$~ $k$-$connectivity$ $rc_k(G)$ of $G$, which is
the minimum integer $j$ for which there exists a $j$-edge-coloring
of $G$ such that for any two distinct vertices $u$ and $v$ of $G$,
there exist at least $k$ internally disjoint $u-v$ rainbow paths.

The concept of rainbow $k$-connectivity has applications in
transferring information of high security in communication networks.
For details we refer to \cite{Chartrand 2} and \cite{Ericksen}.

In \cite{Chartrand 2}, Chartrand et al. studied the rainbow
$k$-connectivity of the complete graph $K_n$ for various pairs $k,
n$ of integers. It was shown in \cite{Chartrand 2} that for every
integer $k\geq 2$, there exists an integer $f(k)$ such that
$rc_k(K_n)=2$ for every integer $n\geq f(k)$. In \cite{Li-Sun 2}, We
improved the upper bound of $f(k)$ from $(k+1)^2$ to
$ck^{\frac{3}{2}}+C$ (here $0<c<1$ and $C=o(k^{\frac{3}{2}})$),
i.e., from $O(k^2)$ to $O(k^{\frac{3}{2}})$. Chartrand et al. in
\cite{Chartrand 2} also investigated the rainbow $k$-connectivity of
$r$-regular complete bipartite graphs for some pairs $k, r$ of
integers with $2\leq k\leq r$, and they showed that for every
integer $k\geq 2$, there exists an integer $r$ such that
$rc_k(K_{r,r})=3$. However, they could not show a similar result as
for complete graphs, and therefore they left an open question: For
every integer $k\geq 2$, determine an integer (function) $g(k)$, for
which $rc_k(K_{r,r})=3$ for every integer $r\geq g(k)$, that is, the
rainbow $k$-connectivity of the complete bipartite graph $K_{r,r}$
is essentially 3. This short note is to solve this question by
showing that $rc_k(K_{r,r})=3$ for every integer $r\geq
2k\lceil\frac{k}{2}\rceil$. We use a method similar to but more
complicated than the proof of Theorem \ref{thm1} in \cite{Chartrand
2}. For notation and terminology not defined here, we refer to
\cite{Bondy}.

\section{Main Result}

In \cite{Chartrand 2}, the authors derived the following results:
\begin{pro}(\cite{Chartrand 2})\label{pro1} For each integer $r\geq 2$,
\[
 rc_2(K_{r,r})=\left\{
   \begin{array}{ll}
     4 &\mbox {if~$r=2$ }\\
     3 &\mbox {if~$r\geq 3$.}
   \end{array}
   \right.
\]\qed
\end{pro}

\begin{pro}(\cite{Chartrand 2})\label{pro2} For each integer $r\geq 3$,
$rc_3(K_{r,r})=3$.\qed
\end{pro}

\begin{thm}(\cite{Chartrand 2})\label{thm1} For every integer $k\geq
2$, there exists an integer $r$ such that $rc_k(K_{r,r})=3$.\qed
\end{thm}

The authors of \cite{Chartrand 2} showed that
$r=2k\lceil\frac{k}{2}\rceil$ is a desired integer for Theorem
\ref{thm1}. We will show, in fact, that $rc_k(K_{r,r})=3$ for every
integer $r\geq 2k\lceil\frac{k}{2}\rceil$, using a method similar to
but more complicated than the proof of Theorem \ref{thm1} in
\cite{Chartrand 2}.

\begin{thm}\label{thm2}For every integer $k\geq 2$, there exists an integer
$g(k)$ such that $rc_k(K_{r,r})=3$ for any $r\geq g(k)$.
\end{thm}
\begin{pf}Let $g(k)=2k\lceil\frac{k}{2}\rceil$. We will show that
$rc_k(K_{r,r})=3$ for every $k\geq 2$, where $r\geq
2k\lceil\frac{k}{2}\rceil$ is an integer. By Propositions \ref{pro1}
and \ref{pro2}, we know that the conclusion holds for $k=2, 3$. So
we assume $k\geq 4$.

We first assume that $k$ is even. Then, $g(k)=2k\cdot{\frac{k}{2}}$.
Since $r\geq g(k)$, then $r=k_1\cdot(2k)+r_1$, where $k_1\geq
\frac{k}{2}, 1\leq r_1\leq 2k-1$. Let the bipartite sets of
$G=K_{r,r}=K_{k_1\cdot(2k)+r_1,k_1\cdot(2k)+r_1}$ be $U$ and $W$.
Let $U'$, $W'$ be the set of first $k_1\cdot(2k)$ vertices of $U$,
$W$, respectively. $U\setminus U'=\{u_1,\ldots,u_{r_1}\}$ and
$W\setminus W'=\{w_1,\ldots,w_{r_1}\}$. Suppose that $$U'=U_1'\cup
\ldots \cup U_{2k}',W'=W_1'\cup \ldots \cup W_{2k}',$$ where
$U_i'=\{u_{i,1},\ldots,u_{i,k_1}\}$ and
$W_i'=\{w_{j,1},\ldots,w_{j,k_1}\}$ for $1\leq i, j\leq 2k$. Let
$G'$ be an induced subgraph of $G$ with bipartite sets $U'$ and
$W'$. Suppose that
$$U=U_1\cup \ldots \cup U_{2k},W=W_1\cup \ldots \cup W_{2k},$$ where
$U_i=U_i'\cup \{u_i\}$, $W_j=W_j'\cup \{w_j\}$ for $1\leq i, j\leq
r_1$ and $U_i=U_i'$, $W_j=W_j'$ for $r_1+1\leq i, j\leq 2k$.

We now give $G$ a $3$-edge coloring as follows: Let $G_1'$ be the
spanning subgraph of $G'$ such that $E(G_1')=\{u_{i,p}w_{j,p}:1\leq
i, j\leq 2k, 1\leq p\leq k_1,i$ and $j$ are of the same parity$\}$.
Let $G_1 $ be the spanning subgraph of $G$ such that
$E(G_1)=E(G_1')\cup \{u_iw_j:1\leq i,j\leq r_1, i$ and $j$ are of
the same parity$\}$. Let $G_2$ be the spanning of subgraph of $G$
such that $$G_2=H_1\cup\ldots\cup H_{2k},$$ where $H_1$ has
bipartite sets $U_1$ and $W_{2k}$, $H_i \ (2\leq i\leq 2k)$ has
bipartite sets $U_i$ and $W_{i-1}$. So,
$H_i=K_{m,n}(\{m,n\}=\{k_1,k_1+1\})$.
\begin{figure}[!hbpt]
\begin{center}
\includegraphics[scale=1.000000]{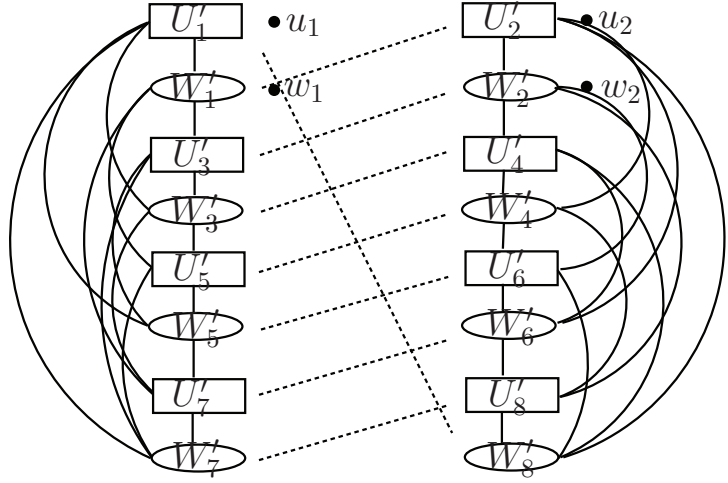}
\end{center}
\caption{The figure for the case $r=18$, $k=4$,
$r_1=2$.}\label{figure1}
\end{figure}
See Figure \ref{figure1} for the case $r=18$, $k=4$, $r_1=2$.
Finally, let
$$G_3=G-(E(G_1)\cup E(G_2)).$$ Assign each edge of $G_i(1\leq i\leq
3)$ the color $i$.

Next we will show that the above edge-coloring is a $k$-rainbow
coloring, that is, there are at least $k$ internally disjoint
rainbow paths connecting any two distinct vertices $u,v$ of $G$. We
will consider the following two cases:

\textbf{Case 1.} $u \in V(G')$. Without loss of generality, let
$u=u_{1,1}$.

\textbf{Subcase 1.1.} $u$ and $v$ belong to the same bipartite set
of $G$.

\textbf{Subsubcase 1.1.1.} $v\in U_1$. Then $G$ contains the $k$
internally disjoint $u_{1,1}-v$ rainbow paths $u_{1,1},w_{i,1},v$
where $1\leq i\leq 2k-1$ and $i$ is odd.

\textbf{Subsubcase 1.1.2.} $v \in U_i$, $3\leq i\leq 2k-1$, and $i$
is odd, say $v\in U_3$. Then $G$ contains the $2k_1\geq k$
internally disjoint $u-v$ rainbow paths $u_{1,1},w_{2,j},v$ and
$u_{1,1},w_{2k,j},v$, where $1\leq j\leq k_1$.

\textbf{Subsubcase 1.1.3.} $v \in U_i$, $2\leq i\leq 2k$, and $i$ is
even, say $v\in U_2$. Then $G$ contains the $2k_1\geq k$ internally
disjoint $u-v$ rainbow paths $u_{1,1},w_{1,j},v$ and
$u_{1,1},w_{2k,j},v$, where $1\leq j\leq k_1$.

\textbf{Subcase 1.2.} $u$ and $v$ belong to different bipartite
sets, and so $v\in W$.

\textbf{Subsubcase 1.2.1.} $v\in W_i$, where $1\leq i\leq 2k-1$ and
$i$ is odd, say $v\in W_1$. Then $G$ contains the $2k_1\geq k$
internally disjoint $u-v$ rainbow paths $u_{1,1},w_{2,j},u_{2,j},v$
and $u_{1,1},w_{2k,j},u_{2k,j},v$, where $1\leq j\leq k_1$.

\textbf{Subsubcase 1.2.2.} $v\in W_i$, where $2\leq i\leq 2k$ and
$i$ is even, say $v\in W_2$. If $v \in W_2'$, without loss of
generality, let $v=w_{2,1}$, then $G$ contains the $u_{1,1}-v$ path
$u_{1,1},v$ together with the $u_{1,1}-v$ rainbow paths
$u_{1,1},w_{3,j},u_{3,j},v$; $u_{1,1},w_{3,1},u_{4,j},v$ and
$u_{1,1},w_{2k,j},u_{2k,j},v$, where $2\leq j\leq k_1$. The cases
for $v=w_2$ and $v\in W_{2k}$ are similar.

\textbf{Case 2.} $u \in V(G)\setminus V(G')$, that is, $u\in
\{u_1,\ldots,u_{r_1};w_1,\ldots,w_{r_1}\}$. Without loss of
generality, let $u=u_1$. By Case 1, we only need to show that there
are at least $k$ internally disjoint rainbow paths connecting $u$
and $v$ for every $v\in V(G)\setminus V(G')$.

\textbf{Subcase 2.1.} $u$ and $v$ belong to the same bipartite set
of $G$.

\textbf{Subsubcase 2.1.1.} $v=u_i$, $3\leq i\leq 2k-1$ and $i$ is
odd, say $v=u_3$. Then $G$ contains the $2k_1\geq k$ internally
disjoint $u-v$ rainbow paths $u_1,w_{2,j},u_3$ and
$u_1,w_{2k,j},u_3$, where $1\leq j\leq k_1$.

\textbf{Subsubcase 2.1.2.} $v=u_i$, $2\leq i\leq 2k$ and $i$ is
even, say $v=u_2$. Then $G$ contains the $2k_1\geq k$ internally
disjoint $u-v$ rainbow paths $u_1,w_{1,j},u_2$ and
$u_1,w_{2k,j},u_2$, where $1\leq j\leq k_1$.

\textbf{Subcase 2.2.} $u$ and $v$ belong to different bipartite sets
of $G$.

\textbf{Subsubcase 2.2.1.} $v=w_i$, $1\leq i\leq 2k-1$ and $i$ is
odd, say $v=w_1$. Then $G$ contains the $2k_1\geq k$ internally
disjoint $u-v$ rainbow paths $u_1,w_{2,j},u_{2,j},w_1$ and
$u_1,w_{2k,j},u_{2k,j},w_1$ where $1\leq j\leq k_1$.

\textbf{Subsubcase 2.2.2.} $v=w_i$, $2\leq i\leq 2k$ and $i$ is
even, say $v=w_2$. Then $G$ contains the $2k_1\geq k$ internally
disjoint $u-v$ rainbow paths $u_1,w_{1,j},u_{3,j},w_2$ and
$u_1,w_{2k,j},u_{2k,j},w_2$ where $1\leq j\leq k_1$.

So the conclusion holds for the case that $k$ is even.

Next we assume that $k$ is odd. Then $g(k)=2k\cdot{\frac{k+1}{2}}$.
Since $r\geq g(k)$, then $r=k_2\cdot(2k)+r_2$, where $k_2\geq
\frac{k+1}{2}, 1\leq r_2\leq 2k-1$. Then with a similar argument to
the case that $k$ is even, we can show that the conclusion also
holds when $k$ is odd. \qed
\end{pf} \\[3mm]
{\bf Remark 2.5.} In \cite{Li-Sun 2} we showed that for every pair
of integers $k \geq 2$ and $r \geq 1$, there is an integer $f(k, r)$
such that if $\ell \geq f(k, r)$, then the rainbow $k$-connectivity
of an $r$-regular complete $\ell$-partite graph is $2$, where
$r$-regular means that every partite set has the same number $r$ of
elements. That is, for sufficiently many number $\ell$ of partite
sets, the rainbow $k$-connectivity of an $r$-regular complete
$\ell$-partite graph is $2$. Theorem \ref{thm2} of this note implies
that for sufficiently large size $r$ of every partite set, the
rainbow $k$-connectivity of an $r$-regular complete $\ell$-partite
graph is at most $3$. So, an interesting question is to think about
the question of determining some bounds on $k,r,\ell$ that tell us
the rainbow $k$-connectivity of an $r$-regular complete
$\ell$-partite graph is $2$ or $3$.


\begin{thebibliography}{99}

\bibitem{Bondy} J.A. Bondy, U.S.R. Murty, {\it Graph Theory}, GTM 244, Springer, 2008.

\bibitem{Chartrand 2} G. Chartrand, G.L. Johns, K.A. McKeon, P. Zhang,
 The rainbow connectivity of a graph, {\it Networks}, 54(2009), 75-81.

\bibitem{Ericksen} A. B. Ericksen, A matter of security, {\it Graduating Engineer
\& Computer Careers} (2007), 24-28.

\bibitem{Li-Sun 2} X. Li, Y. Sun, The rainbow $k$-connectivity
of two classes of graphs, arXiv:0906.3946 [math.CO].





\end{thebibliography}
\end{document}